\documentclass{amsart}
\usepackage{stmaryrd}
\usepackage{amssymb}
\usepackage{amsmath}
\usepackage{amsthm}
\usepackage{amsbsy}
\usepackage[dvips]{graphics}
\usepackage{amscd}
\usepackage{bm}

\theoremstyle{plain}
\newtheorem{theorem}{Theorem}[section]

\theoremstyle{definition}

\theoremstyle{remark}

\newtheorem*{acknowledgements}{Acknowledgements}

\newcommand{\C}{{\mathbb{C}}}

\begin{document} 

\title[Negative isotropic curvature]{A note on Negative isotropic curvature}

\author{Harish Seshadri}
\email{harish@isibang.ac.in}
\address{Stat-Math Unit,\\
Indian Statistical Institute,\\
Bangalore, India}

\date{\today}

\begin{abstract}
We prove that any smooth orientable closed four-manifold admits a
Riemannian metric with negative isotropic curvature in the sense of
Micallef and Moore.
\end{abstract}

\maketitle
\section{Introduction}
An interesting notion of curvature, that of {\it curvature on
isotropic 2-planes}, was introduced by Micallef and Moore in
~\cite{mm}. The definition is as follows. Let $(M,g)$ be a Riemannian
manifold of dimension at least $4$. Consider the complexification $TM
\otimes _{\Bbb R} {\Bbb C}$. The inner products on $TM$ and 
$\wedge ^2 TM$ extend to complex bilinear maps on
$TM \otimes \C$ and $\wedge ^2 TM \otimes \C$ and we denote these extensions
by $\langle, \rangle$ in both cases. A subspace $W \subset T_pM
\otimes {\Bbb C}$ is said to be {\it isotropic} if $\langle w,w
\rangle=0$ for all $w \in W$. Now let $R: \wedge ^2 TM
\otimes \C \rightarrow \wedge ^2 TM \otimes \C$ denote
the complex linear extension of the curvature
operator of $g$. The metric $g$ is said to have {\it
positive} (resp. negative) {\it isotropic curvature} if $\langle R(v
\wedge w), \overline v \wedge \overline w \rangle >0$ (resp.$<0$)
whenever $\{ v,w \}$ span an isotropic $2$-plane.

Considerations of the second variation of area of minimal surfaces
lead naturally to isotropic curvature. Examples of manifolds with
isotropic curvature of a fixed sign are pointwise quarter-pinched
manifolds and manifolds with positive or negative curvature operators.

The relationship between the sign of isotropic curvature and the
``classical curvatures'' is somewhat mysterious. It is known
(cf. ~\cite{mw}) that positive isotropic curvature does not imply
positive Ricci curvature but does imply positive scalar curvature (a
similar statement is true for negative isotropic curvature). As in the
case of positive scalar curvature, one knows that if two manifolds
have positive isotropic curvature then their connected sum admits a
metric with positive isotropic curvature. However, unlike the case of
positive scalar curvature, surgeries in lower codimension need not
preserve positive isotropic curvature.

Pursuing the analogy between isotropic and scalar curvatures, one is
tempted to ask whether every smooth closed manifold admits a metric
with negative isotropic curvature. In this note we show that this is
the case in dimension $4$.

\begin{theorem}\label{mai}
Any smooth orientable closed four-manifold admits a $C^{ \infty }$
metric with negative isotropic curvature.
\end{theorem}

{\it Remarks}:

(i) A K\"ahler metric can never have negative isotropic curvature.
In fact, if $v$ and $w$ are such that $(v \wedge w)^{1,1}=0$, then
$\langle R(v \wedge w), \overline v \wedge \overline w \rangle =0$.
There exist $v$ and $w$ such that $(v \wedge w)^{1,1}=0$ and $\C \{ v,
w \}$ is an isotropic 2-plane ~\cite{sou}.

(ii) In the definition of isotropic curvature, if one demands that
$\langle R(v \wedge w), v \wedge w \rangle <0$ for {\it any} $v \wedge
w \in \wedge ^2 TM \otimes {\Bbb C}$, then one says, following
Sampson ~\cite{sam}, that $g$ has negative Hermitian sectional curvature. This
curvature occurs in the study of harmonic maps of K\"ahler manifolds
and clearly implies negative sectional curvature. Hence it has strong
topological implications while Theorem \ref{mai} states that negative
isotropic curvature, in dimension 4 at least, has {\it no} topological
implications. \\

The proof of Theorem \ref{mai} essentially hinges on the fact that in
dimension four, the negative isotropic curvature condition can be
rephrased as a variational problem. The functional involved will be
the integral of a scalar curvature-like function $\sigma$ and its
domain will be the space of metrics. Once this is noted the proof
involves showing that the functional becomes negative for a suitable
choice of metric $g$. Then standard techniques in conformal geometry
(as in ~\cite{gur}) can be applied to show that there exists a
conformal deformation of $g$ for which $\sigma $ is negative. The
geometric part of the proof is the construction of $g$. We do this by
gluing in a product of a circle with a hyperbolic 3-manifold of finite
(but large) volume. This hyperbolic 3-manifold will be a knot
complement in $S^3$ and the gluing will be done without changing the
smooth topology of $M$.

Note that the proof below only gives a $C^{2, \alpha}$ metric with negative
isotropic curvature. The density of $C^ \infty$ metrics in the space of 
$C^{2,\alpha}$ metrics (with the $C^{2, \alpha}$ norm) will then give us
a $C^ {\infty}$ metric with negative isotropic curvature.

\begin{acknowledgements}
I would like to thank the referee for pointing out various errors in the
original draft.
\end{acknowledgements}

\section{Proof}   
\begin{proof}
For the sake of exposition we divide the proof into two steps.\\

STEP 1: Let $M$ be a closed smooth four-manifold. A simple calculation, as in
\cite{mm}, shows that a metric $g$ on $M$ has negative isotropic
curvature if and only if 
\begin{equation}
{\frac {s}{6}} I-W < 0,  \label{one}
\end{equation}
where $s$ denotes the scalar
curvature, $I$ is the identity operator and $ W$ the trace-zero
symmetric operator given by Weyl curvature on $\wedge^2TM$. To prove 
(\ref{one}) it is enough to show that
$$ \sigma_g:= {\frac {s}{6}} + \vert W \vert <  0,$$ where $\vert W
\vert$ denotes the pointwise norm of $W$ (regarded as an element of 
$ \otimes ^4 T^*M$). This is because $W$ is a symmetric operator and
hence $ \langle W( \alpha ), \alpha \rangle \ge - \vert W \vert \vert
\alpha \vert ^2$.

The main property of $ \sigma _g$ is that under conformal changes it
transforms like scalar curvature. More precisely, if $g \rightarrow
\tilde g=u^{\frac{4}{n-2}} g$, the modified scalar curvature
transforms (in what follows, the Laplacian $\triangle =-d^ \ast d$) by
\begin{equation}
\sigma _g \rightarrow \sigma _{\tilde g}  = u^{\frac{-4}{n-2}} 
\sigma _g -4{\frac{n-1}{n-2}} u^{- {\frac{n+2}{n -2}}} \triangle u. \label{tra}
\end{equation}
This follows immediately from the corresponding transformation law for
the usual scalar curvature $s$ cf. ~\cite{bes} and the fact that
$\vert W_{f^2g} \vert _{f^2g}=f^{-2} \vert W_g \vert_g$ for any smooth 
positive function $f$.

Consider the functional on the space of $C^2$ metrics on $M$ defined by
$$F(g)= \int _M \sigma _g dV_g, $$ where $dV_g$ is the volume form of
the metric $g$. A standard argument, which we outline below, shows
that if $F(g) <0$ for some $ g$, then there is a $C^{2, \alpha}$
metric $\tilde g$ in the conformal class of $ g$ with ${\tilde
\sigma_g} <0$.

Define the operator $L$ on $C^\infty(M)$ by
$$L=-4{\frac {(n-1)}{(n-2)}}\triangle + \sigma _g.$$
Let $<,>$ denotes the $L^2$ inner product on $C^\infty(M)$.
and let
$$\lambda= inf _{{f \in W^{1,2}} \atop {\Vert f \Vert _2=1}}<Lf,f>$$
and $u$ be the corresponding eigenfunction. We note that since
$\sigma _g$ is, in general, Lipschitz continuous
but not smooth (at the zero locus of $\vert W \vert$), the best regularity we
can obtain for $u$ is that $u \in C^{2, \alpha}$ for any $0< \alpha <1$. This
is sufficient for our purposes. By the minimum principle $u>0$ and
by definition, $u$ satisfies 
\begin{equation}
 L(u)= \lambda u. \label{pos}
\end{equation} 
Now suppose that $g$ satisfies $F(g) <0$. Then
$< L(1),1> =F(g)<0$. Hence $\lambda <0$. 

Consider the
metric $\tilde g= u^{\frac{4}{n-2}} \tilde g$. By (\ref{tra}) and
(\ref{pos}), we see that
$$ \sigma _g= u^{-{\frac {(n-2)}{4(n-1)}}-1} L(u) <0.$$ 

Hence we are reduced to constructing a metric $g$ with $F(g)<0$. 
This is done in the next step. \\

STEP 2: Let $K$ be a knot in $S^3$ such that
$N:=S^3 - K$ admits a complete hyperbolic metric (constant curvature
$-1$) $g_H$ of finite volume. Such a knot exists by the work of Riley
\cite{ril}. Let $ \overline N =S^3 - \nu (K) \subset N$, where $ \nu
(k)$ is an open tubular neighbourhood of $K$. Fix an embedding of
$\overline N \times S^1$ in $M$ (which exists since $S^3-K$ can be
embedded in a 3-ball).

Briefly, we will construct metrics, depending on a parameter $c$, on $
\overline N \times S^1$ which are isometric to a fixed standard metric
near $\partial ( \overline N \times S^1)$. Then these metrics can be
extended to all of $M$ such that the extensions agree on $M-(\overline
N \times S^1$). For large values of $c$ we get metrics for which $F
<0$.
 
The actual construction is as follows: $N$ has one cusp end and this end
is isometric to $(E:=[0, \infty) \times T^2, g_H:=dt^2 +e^{-2t}g_e)$,
where $T^2$ is the 2-torus and $g_e$ a flat metric on it. Now consider
the metrics $c^2g_H$, where $c \in {\Bbb R}^+$ on $N$. We note two
features of these metrics.

(i) The sectional curvatures of these metrics go to zero as $c \rightarrow
\infty$. Also, by choosing $a(c)$ appropriately (in fact
$a(c)=c \ log({\frac {c}{2}})$ will work, see $(\ref{bd})$ below), we can 
ensure that the ``cross-sectional'' torii are approximately 
of fixed size on the portion $[a(c),a(c)+1] \times T^2$ of the end.

(ii) For any $\alpha \ge 0$, let $N_{ \alpha } :=N- \bigl ( [ \alpha, \infty)
    \times T^2 \bigr)$. Then 
\begin{equation}\label{sca}
   \int_{ N_{0}}  s_{c^2g_H}dV_{c^2g_H} =
    -6c Vol(N_0, g_H) \rightarrow - \infty \ {\rm as} \ 
c \rightarrow  \infty.
\end{equation}

Now the scaled metric $c^2g_H$ on the end can be written as
$$c^2g_H:=dt^2 + c^2e^{-2 {\frac {t}{c}}}g_e \ \ {\rm with} \ \ t \in 
[0, \infty).$$
We deform $c^2g_H$ to the product (flat) metric
$dt^2 +g_e$ in the interval $[a(c),a(c)+1]$, $a(c)$ to be
specified, as follows: let $ \phi : {\Bbb R} \rightarrow [0,1]$ be a
smooth decreasing function such that
\begin{align} \notag
\phi (t) & = 1 \ \ {\rm for} \ \ t \le 0 \\ \notag
         & = 0 \ \ {\rm for} \ \ t \ge {\frac {1}{2}}. \\ \notag
\end{align}  
Further, let $$f_c(t)= \phi \bigl (t-a(c) \bigr )ce^{- {\frac {t}{c}}}
+1- \phi \bigl (t-a(c) \bigr ).$$ Then the metric
\begin{equation} \label{me}
k_c =dt^2 + f_c(t)^2g_e, \ \ t \in [0, \infty) 
\end{equation}
deforms $c^2g_H$ to $dt^2+g_e$ over the interval $[a(c),a(c)+1]$. 

For $c > {\frac {1}{log(2)}}$, let $a(c)=c \ log (c)$. If $t \in
[a(c),a(c)+1]$, then ${\frac {1}{2}} \ < \ ce^{-{\frac {t}{c}}} \ < \
1$. From this an easy computation shows that the $C^2$
distance between the metrics $k_c$ (defined in $(\ref{me})$) and $dt^2
+ g_e$ on $[a(c), a(c)+1] \times T^2$ is bounded,
\begin{equation} \label{bd}
\vert k_c-(dt^2+g_e) \vert_{C^2} \le A,
\end{equation}
with $A$ independent of $c$.

The metric $g_c$ that we are 
interested (on $M$) is now given by: 
\begin{align} \notag
g_c & = c^2g_H + d \theta ^2 \ \ {\rm on} \ \ N_0 \times S^1 \\ \notag
    & = k_c + d \theta ^2 \ \ {\rm on} \ \ [0, a(c)+1 ] \times T^2
    \times S^1 \\ \notag
    & = \overline g \ \ {\rm on} \ \ M - (N_{a(c)+1} \times S^1), \notag
\end{align}    
where $\overline g$ is a fixed extension of the standard product
metric $dt^2 + g_e + d \theta ^2$ on a tubular neighbourhood of $ T^2
\times S^1 \cong \partial \Bigl (M - (\overline N \times S^1) \Bigr )
\cong \partial \Bigl (M - (N_{a(c)+1} \times S^1) \Bigr )$, where $
\cong$ denotes diffeomorphism, to all of $M - ( \overline N \times
S^1)$. This extension can obviously be chosen to be independent of
$c$.

Now we claim that $F(g_c)<0$ for large $c$. To see this, note first
that the scalar curvature integral $\int_M s_{g_c}dV_{g_c}$ is equal
to
\begin{equation}\label{er}
\int_{ N_{a(c)} \times S^1}
 s_{g_c}dV_{g_c}
+ \int_{[a(c),a(c)+1] \times T^2 \times S^1} s_{g_c}dV_{g_c} +
\int_{M- ( \overline N \times S^1)} s_{g_c}dV_{g_c}, 
\end{equation}  
where we have identified $N_{a(c)+1}$ with $ \overline N$ in the last
integral. Now, as $c \rightarrow \infty$, the first integral goes to
$- \infty$ by (\ref{sca}). The second remains bounded by (\ref{bd})
and the third is independent of $c$. Hence
\begin{equation}\label{fir}
\int_M s_{g_c}dV_{g_c} \rightarrow - \infty \ \ {\rm as} \ \ c 
 \rightarrow \infty.
\end{equation}
Similarly, the Weyl curvature integral $\int_M \vert W_{g_c} \vert dV_{g_c}$
is equal to 
\begin{equation}\label{we}
\int_{ N_{a(c)} \times S^1}
 \vert W_{g_c} \vert dV_{g_c}
+ \int_{[a(c),a(c)+1] \times T^2 \times S^1} \vert W_{g_c} \vert dV_{g_c} +
\int_{M- (\overline N \times S^1)} \vert W_{g_c} \vert dV_{g_c},
\end{equation}
where the norm of $W_{g_c}$ is with respect to $g_c$. Again, for the
same reasons, the second and third integrals will remain bounded as $c
\rightarrow \infty$. To see that the first integral is zero, recall
that the metric on $ { N_{a(c)} \times S^1}$ is the product
metric and the product of a constant sectional curvature metric with
$S^1$ is conformally flat. Hence
\begin{equation}\label{sec}
\int_M \vert W_{g_c} \vert 
dV_{g_c} \ \ {\rm remains \ bounded \ as} \ \ c \rightarrow \infty. 
\end{equation}
Combining $(\ref{fir})$ and $(\ref{sec})$ we see that $F(g_c)<0$ for
large $c$.

\end{proof}

\noindent {\it Remarks}: 

(i) The proof above actually shows that
on any closed $4$-manifold, for any $ \mu >0$, there exists a metric
$g$ with $\mu s_g + \vert W_g \vert <0$.

(ii) In the case of negative scalar curvature, one knows by the work
of Lohkamp \cite{loh}, that the space of such metrics $\mathfrak{M}$
is contractible. One can ask if the same is true for the space of
negative isotropic curvature metrics. Note that in Lohkamp's paper, a
variational approach similar to Step 1 is used to
extend a map of $S^{n}$ into $\mathfrak{M}$ to a map
of $B^{n+1}$ (in order to prove that $\pi_n( \mathfrak{M})=0$).

(iii) The proof hinges on the fact that in dimension 4, there exists
an expression of the form $\mu s_g \mp \vert W_g \vert$, with $\mu >0$
such that if this expression has a certain fixed sign, then the
isotropic curvature also has the same sign.  In order to extend the
proof to dimensions $n \ge 5$, one needs to know if a similar
statement is true for all $n$. However the following example, given in
\cite{mw} shows that such a statement is not true in the positive
case. The product metric on $M= \Sigma \times S^{2m}$, where $\Sigma$
is a compact surface of constant curvature $-1$ and $S^{2m}$ is the
unit sphere of dimension $2m \ge 4$, is conformally flat and has
positive scalar curvature. However, since $H_2(M, {\mathbb R}) \neq 0$
by Theorem 2.1 of ~\cite{mw}, $M$ does not admit a metric of positive
isotropic curvature.

\end{document}